\documentclass{ifacconf}

\counterwithin*{section}{part}

\usepackage[english]{babel}
\usepackage{amssymb}
\usepackage{tikz}
\usepackage{mathtools}
\usepackage{etoolbox}
\usepackage{xifthen}
\usepackage{enumitem}
\usepackage{stmaryrd}

\usepackage[ruled,vlined,algo2e,linesnumbered]{algorithm2e}
\makeatletter
\patchcmd{\algocf@Vline}{\vrule}{\vrule\vspace{-.32em}}{}{}
\makeatother
\SetStartEndCondition{ }{ }{}%
\SetKw{KwTo}{to}\SetKwFor{For}{for}{\string do}{}%
\SetKwIF{If}{ElseIf}{Else}{if}{then}{elif}{else}{}%
\SetKwFor{While}{while}{\string do}{}%

\DeclareMathAlphabet{\pazocal}{OMS}{zplm}{m}{n}
\usepackage [autostyle, english = american]{csquotes} 
\MakeOuterQuote{"}

\usepackage{pgfplots}
\pgfplotsset{compat=newest}
\pgfplotsset{plot coordinates/math parser=false}
\newlength\figureheight
\newlength\figurewidth 

\usepackage{xargs}
\usepackage{soul, color, xcolor}
\usepackage[numbers]{natbib}

\usetikzlibrary{positioning,arrows,petri,calc,decorations.markings,arrows.meta}
\tikzset{
place/.style={circle,thick,minimum size=4mm,draw},
transitionV/.style={rectangle,thick,fill=black,minimum height=6mm,inner xsep=1pt}
}

\definecolor{myblue}{RGB}{0, 101, 202}
\definecolor{mygreen}{RGB}{130, 180, 0}
\definecolor{myred}{RGB}{197, 14, 31}
\definecolor{mypurple}{RGB}{128, 0, 128}
\definecolor{myyellow}{HTML}{F7E2B0}
\definecolor{mygrey}{RGB}{105, 105, 105}

\newcommand{\type}{j}
\newcommand{\types}{\MakeUppercase{\type}}
\newcommand{\machine}{m}
\newcommand{\machines}{\MakeUppercase{\machine}}
\newcommand{\quantity}{k}
\newcommand{\quantities}{\MakeUppercase{\quantity}}
\newcommand{\capacity}{C}
\newcommand{\typeset}{\dint{1,\types}}
\newcommand{\machineset}{\dint{1,\machines}}

\newcommand{\quantityset}{\dint{1,\quantities}}
\newcommand{\batch}{b}
\newcommand{\batches}{\MakeUppercase{\batch}}

\newcommand{\dint}[1]{\left\llbracket#1\right\rrbracket} 

\newcommand{\N}{\mathbb{N}}
\newcommand{\No}{\mathbb{N}_0}
\newcommand{\Z}{\mathbb{Z}}

\newcommand{\R}{\mathbb{R}}
\newcommand{\Rmax}{{\R}_{\normalfont\fontsize{7pt}{11pt}\selectfont\mbox{max}}}
\newcommand{\Rmin}{{\R}_{\normalfont\fontsize{7pt}{11pt}\selectfont\mbox{min}}}
\newcommand{\Rbar}{\overline{\R}}

\newcommand{\floor}[1]{\left\lfloor#1\right\rfloor}
\newcommand{\ceil}[1]{\left\lceil#1\right\rceil}

\renewcommand{\qed}{\hfill\ensuremath{\blacksquare}}

\newcommand{\arcs}{E}
\newcommand{\nodes}{N}

\makeatletter
\newcommand{\splus}{%
  \DOTSB\mathop{\mathpalette\mattos@splus\relax}\slimits@
}
\newcommand\mattos@splus[2]{%
  \vcenter{\hbox{%
    \sbox\z@{$#1\oplus$}%
    \resizebox{!}{0.9\dimexpr\ht\z@+\dp\z@}{\raisebox{\depth}{$\m@th#1\boxplus$}}%
  }}%
  \vphantom{\oplus}%
}
\makeatother

\makeatletter
\newcommand{\stimes}{%
  \DOTSB\mathop{\mathpalette\mattos@stimes\relax}\slimits@
}
\newcommand\mattos@stimes[2]{%
  \vcenter{\hbox{%
    \sbox\z@{$#1\otimes$}%
    \resizebox{!}{0.9\dimexpr\ht\z@+\dp\z@}{\raisebox{\depth}{$\m@th#1\boxtimes$}}%
  }}%
  \vphantom{\otimes}%
}
\makeatother

\usepackage{pict2e}

\makeatletter
\newcommand*{\bigsplus}{\DOTSB\mathop{\mathpalette\big@boxplus\relax}\slimits@}

\newcommand{\big@boxplus}[2]{%
  \vcenter{%
    \m@th\bigbox@thickness{#1}%
    \sbox\z@{$#1\bigoplus$}%
    \dimen@=\ht\z@ \advance\dimen@\dp\z@
    \hbox{%
      \setlength{\unitlength}{\dimen@}%
      \begin{picture}(1,1)
      \polyline(0.1,0.1)(0.9,0.1)(0.9,0.9)(0.1,0.9)(0.1,0.1)(0.5,0.1)
      \polyline(0.5,0.1)(0.5,0.9)
      \polyline(0.1,0.5)(0.9,0.5)
      \end{picture}%
    }%
  }%
}

\newcommand{\bigbox@thickness}[1]{%
  \ifx#1\displaystyle
    \linethickness{0.2ex}%
  \else
    \ifx#1\textstyle
      \linethickness{0.16ex}%
    \else
      \ifx#1\scriptstyle
        \linethickness{0.12ex}%
      \else
        \linethickness{0.1ex}%
      \fi
    \fi
  \fi
}
\makeatother
\newcommand{\wA}{\mathsf{a}}
\newcommand{\wB}{\mathsf{b}}
\newcommand{\wC}{\mathsf{c}}

\newcommand{\graph}{\pazocal{G}}
\newcommand{\nonegset}{\Gamma} 

\newcommand{\svdots}{\raisebox{3pt}{$\scalebox{.75}{\vdots}$}} 
\newcommand{\sddots}{\raisebox{3pt}{$\scalebox{.75}{$\ddots$}$}} 

\begin{document}
 

\begin{frontmatter}

\title{Switched max-plus linear-dual inequalities\\for makespan minimization:\\the case study of an industrial bakery shop\thanksref{footnoteinfo}} 

\author[1]{Davide Zorzenon}
\author[2]{Nataliia Zaiets}
\author[1,3]{J\"{o}rg Raisch}

\thanks[footnoteinfo]{Support from Deutsche Forschungsgemeinschaft (DFG) via grant RA 516/14-1 and under Germany’s Excellence Strategy -- EXC 2002/1 ``Science of Intelligence'' -- project number 390523135 is gratefully acknowledged.}

\address[1]{Control Systems Group, Technische Universit\"at Berlin, Germany (e-mail: [zorzenon,raisch]@control.tu-berlin.de)}
\address[2]{Department of Automation and Robotics Systems, National University of Life and Environmental Sciences of Ukraine, Kyiv (e-mail: z-n@ukr.net)}
\address[3]{Science of Intelligence, Research Cluster of Excellence, Berlin, Germany}

\begin{abstract}
In this paper, an industrial bakery shop is modeled by switched max-plus linear-dual inequalities (SLDIs). 
SLDIs are timed discrete event systems suitable for describing flow shops with time-window constraints and switching operating modes, where each mode corresponds to a job type.
We consider the scheduling problem of minimizing the makespan of the shop, and we show that the application of methods based on the max-plus algebra leads to a faster solution compared to standard techniques.
The results of the paper are general, in the sense that they can be applied to any permutation flow shop with time-window constraints.
\end{abstract}


\end{frontmatter}

\section{Introduction}

Industrial bakeries are complex manufacturing systems characterized by a high diversity of products, high production volumes, and time-sensitive processing stages.
In fact, to ensure high-quality products, operations need to be executed under strict temporal requirements; for instance, excessive duration of the yeast fermentation may cause an undesirable collapse of the dough matrix (overproofing).
Most of the energy consumption in the bakery industry comes from the baking ovens.
Thus, the use of this equipment must be carefully planned to reduce idle time, temperature dispersion and, consequently, higher energy costs derived from re-heating.

In this paper, we consider the problem of minimizing the makespan of an industrial bakery shop to increase its efficiency and reduce energy consumption.
The shop, schematically depicted in Figure~\ref{fi:bakery_shop}, is representative of a typical production line.
Experimental data and a list of the main technological equipment were provided by the industrial bakery LLC "Novi Perspektivi", Rivne region, Ukraine.
To simplify the discussion, we assume that no machines working in parallel are present in the shop and defer the analysis of a more general scenario to future work.
We emphasize, though, that such a general scenario is well within the scope of the methods discussed in the following.

To solve the problem, the system is modeled by means of switched max-plus linear-dual inequalities (SLDIs, see~\cite{zorzenon2022switched}).
SLDIs are the extension of switching max-plus linear systems (\cite{VANDENBOOM20061199}) to the case of systems with time-window constraints.
More precisely, they are discrete event systems whose nondeterministic dynamics is described by a sequence of operating modes.
Each mode corresponds to a system of inequalities that restrict the admissible time of occurrence of events, which represent the start or end of a process. 

The makespan minimization problem in bakery systems has been considered, e.g., in~\cite{hecker2014application,babor2021optimization}.
In contrast to these papers, we do not impose the \textit{no-wait} requirement, which forces a process to start as soon as the previous one has ended; this constraint has the advantage of simplifying the solution of the scheduling task, but it may be too restrictive in practice, leading to a suboptimal makespan.
We mention that the description of a bakery plant by means of (non-switched) linear-dual inequalities has already been presented, e.g., in~\cite{declerck2021critical}.
The main improvement of the SLDI model of the present paper is the ability to describe dynamics of different types of products; this is necessary, as different product types require different baking times, for example.

After some mathematical preliminaries (Section~\ref{se:2}), in Section~\ref{se:3} we define and interpret SLDIs in the context of manufacturing systems; based on that, we present a closed formula for the makespan of permutation flow shops with time-window constraints (defined in Section~\ref{su:makespan_permutation}). 
The formula uses only basic max-plus operations, and its implementation can lead to important computational savings with respect to other approaches.
This is demonstrated in Section~\ref{se:4}, where it is applied to the makespan minimization of the bakery shop and compared to standard algorithms from linear programming and graph theory.
Section~\ref{se:conclusions} provides conclusions and suggestions for future work.

\textbf{Notation.}
We denote sets $\R\cup\{-\infty\}$, $\R\cup\{+\infty\}$, $\R\cup\{\pm\infty\}$ respectively by $\Rmax$, $\Rmin$, and $\Rbar$.
The sets of nonnegative and positive integers are denoted, respectively, by $\No$ and $\N$.
For a matrix $A\in\Rbar^{n\times n}$, $A^\sharp$ indicates $-A^\top$.
Given $a,b\in\No$ with $b\geq a$, $\dint{a,b}$ denotes the discrete interval $\{a,a+1,a+2,\ldots,b\}$.
Moreover, given $x\in\R$, $\ceil{x}\coloneqq \min \{y\in\Z\ | \ y\geq x\}$ and $\floor{x}\coloneqq \max\{y\in\Z\ |\ y\leq x\}$.
\section{Preliminaries}\label{se:2}

In this section we recall some preliminaries on max-plus algebra and precedence graphs; for a more detailed overview of these topics we refer to~\cite{baccelli1992synchronization,butkovivc2010max,hardouin2018control}.

\subsection{Max-plus algebra}

The max-plus algebra is the mathematical framework consisting of real numbers extended with $-\infty$ and $+\infty$, and operations $\oplus$ (addition), $\otimes$ (multiplication), $\splus$ (dual addition), and $\stimes$ (dual multiplication) defined as follows: for all $a,b\in\Rbar$,
\[
\begin{array}{lcl}
    a \oplus b = \max(a,b), && a\otimes b = 
    \begin{dcases}
        a+b & \mbox{if } a,b\neq -\infty,\\
        -\infty & \mbox{otherwise,}
    \end{dcases}\\
    a \splus b = \min(a,b), && a\stimes b = 
    \begin{dcases}
        a+b & \mbox{if } a,b\neq +\infty,\\
        +\infty & \mbox{otherwise.}
    \end{dcases}
\end{array}
\]

These operations can be naturally extended to matrices; given $A,B\in\Rbar^{m\times n}$, $C\in\Rbar^{n\times p}$, for all $i\in\dint{1,m}$, $j\in\dint{1,n}$, $h\in\dint{1,p}$,
\[
    \begin{array}{lcl}
        (A\oplus B)_{ij} = A_{ij}\oplus B_{ij}, && (A\otimes C)_{ih} = \bigoplus_{k=1}^n A_{ik}\otimes C_{kh},\\
        (A\splus B)_{ij} = A_{ij}\splus B_{ij}, && (A\stimes C)_{ih} = \bigsplus_{k=1}^n A_{ik}\stimes C_{kh}.
    \end{array}
\]
Matrices $\pazocal{E}$, $\pazocal{T}$, and $E_\otimes$ are, respectively, the neutral element for $\oplus$, $\splus$, and $\otimes$, i.e., $\pazocal{E}_{ij} = -\infty$ and $\pazocal{T}_{ij}=+\infty$ for all $i,j$, and $(E_\otimes)_{ij} = 0$ if $i=j$ and $(E_\otimes)_{ij} = -\infty$ if $i\neq j$.
The $r$th power of square matrix $A$ is defined recursively by $A^0 = E_\otimes$ and $A^r = A^{r-1} \otimes A$. 
The Kleene star of $A$ is $A^* = \bigoplus_{i \in\No} A^i$.
The partial order relation $\preceq$ between two matrices of the same dimension is induced by $\oplus$ as: $A\preceq B\ \Leftrightarrow\ A\oplus B = B$; hence, $A\preceq B$ is equivalent to $A_{ij}\leq B_{ij}$ for all $i,j$.

\subsection{Precedence graphs}

Precedence graphs are commonly used in scheduling theory to graphically represent the time relations between processes in various types of manufacturing environments (\cite{pinedo2016scheduling}); here we recall their connection with the max-plus algebra.
Given a matrix $A\in\Rmax^{n\times n}$, the precedence graph corresponding to $A$ is the pair $\graph(A)=(\nodes,\arcs)$, where $\nodes=\dint{1,n}$ is the set of nodes, and $\arcs\subseteq \nodes\times\nodes$ is the set of weighted arcs, defined such that there is an arc $(j,i)\in\arcs$ with weight $A_{ij}$ if and only if $A_{ij}\neq-\infty$.

A path on $\graph(A)$ is a sequence of nodes $\rho=(i_1,i_2,\ldots,i_{r+1})$ such that $(i_j,i_{j+1})\in\arcs$ for all $j\in\dint{1,r}$; the number $|\rho|=r$ is the length of $\rho$.
A circuit is a path in which the first and last node coincide, i.e., $i_1=i_{r+1}$.
The weight $w_\rho$ of path $\rho$ is the sum (in standard algebra) of the weights of the arcs composing it; in the max-plus algebra, this quantity can be computed as $w_\rho=\bigotimes_{j=1}^{|\rho|}A_{i_{j+1}i_j}$.
Generalizing the latter formula, we get the following graphical interpretation of the Kleene star of a matrix: $(A^*)_{ij}$ is equal to the maximum weight of all the paths in $\graph(A)$ from node $j$ to node $i$.
Matrix $A^*$ belongs to $\Rmax^{n\times n}$ if and only if there are no circuits with positive weight in $\graph(A)$, otherwise at least one element in the diagonal of $A^*$ is $+\infty$.
We denote by $\nonegset$ the set of precedence graphs with no circuits with positive weight.
Property "$\graph(A)\in\nonegset$" can be checked in time $\pazocal{O}(n^3)$, and, if the property holds, $A^*$ can be computed in the same time complexity.


\section{Scheduling problems in the max-plus algebra}\label{se:3}

In this section, we characterize the dynamics of manufacturing systems with time-window constraints by means of switched max-plus linear-dual inequalities (SLDIs).
After that, we show how to solve makespan minimization problems in the max-plus algebra.

\subsection{Switched max-plus linear-dual inequalities}

Let $\Sigma=\{\wA_1,\ldots,\wA_m\}$ be a finite set of \textit{modes}, and let $v = (v_1,v_2,\ldots ,v_{\quantities})$ be a finite ordered sequence of modes\footnote{In~\cite{zorzenon2022switched}, SLDIs were defined for the case of infinite sequences of modes; in this paper, only finite sequences will be considered.} $v_k\in\Sigma$, with $k\in\dint{1,\quantities}\subset\N$.
SLDIs are systems of inequalities in the variables $x(1),\ldots,x(\quantities)\in\R^n$ of the following form:
\begin{equation}\label{eq:SLDI}
    \begin{array}{rc}
    \forall k\in\dint{1,\quantities}: & A^0_{v_k} \otimes x(k) \preceq x(k) \preceq B^0_{v_k} \stimes x(k),\\
    \forall k\in\dint{1,\quantities-1}: & A^1_{v_k} \otimes x(k) \preceq x(k+1) \preceq B^1_{v_k} \stimes x(k),
    \end{array}
\end{equation}
where $A^0_{\wA_i},A^1_{\wA_i}\in\Rmax^{n\times n}$ and $B^0_{\wA_i},B^1_{\wA_i}\in\Rmin^{n\times n}$ for all $i\in\dint{1,m}$.

A possible interpretation of~\eqref{eq:SLDI} is as follows.
Let $x_i(k)$ represent the time of the $k$th occurrence of event $i\in\dint{1,n}$; event $i$ could for instance be the start or the end of a process in a shop.
The temporal difference between the occurrence of two events may be subject to job-dependent window constraints; typically, lower bounds describe the time required for completing a process, whereas upper bounds represent quality specifications.
If we associate each job type to a mode $\wA_i$, then sequence $v$ describes the entrance order of jobs into the shop, matrices $A^0_{v_k}$, $A^1_{v_k}$ contain the temporal lower bounds, and matrices $B^0_{v_k},B^1_{v_k}$ the upper bounds referred to job of type $v_k$.
The first pair of inequalities in~\eqref{eq:SLDI} specifies that $x_i(k)$ must be between $\max_j \big((A^0_{v_k})_{ij} + x_j(k)\big)$ and $\min_j \big((B^0_{v_k})_{ij} + x_j(k)\big)$; the second pair of inequalities relates $x_i(k+1)$ with $x_j(k)$ (i.e., the $(k+1)$st occurrence of event $i$ with the $k$th occurrence of event $j$), for all $j$, in a similar way.
Since, for all $a,b,c\in\Rbar$, $\max(a,b) \leq c$ is equivalent to $a \leq c$ and $b \leq c$, and $a \leq \min(b,c)$ is equivalent to $a \leq b$ and $a\leq c$, we can rewrite~\eqref{eq:SLDI} as: for all $i,j\in\dint{1,n}$, $k\in\dint{1,\quantities}$,
\begin{equation*}
    (A^0_{v_k})_{ij} \leq x_i(k) - x_j(k) \leq (B^0_{v_k})_{ij},
\end{equation*}
and, for all $i,j\in\dint{1,n}$, $k\in\dint{1,\quantities-1}$,
\begin{equation*}
    (A^1_{v_k})_{ij} \leq x_i(k+1) - x_j(k) \leq (B^1_{v_k})_{ij},
\end{equation*}
which is a linear system of inequalities.
Any sequence $x(1),\ldots,x(\quantities)$ that satisfies the latter inequalities represents a valid (or consistent) trajectory of jobs in the shop, in the sense that lower and upper bound constraints are always satisfied.

\subsection{Makespan minimization in permutation flow shops}\label{su:makespan_permutation}

A typical objective considered in scheduling theory is the minimization of the makespan in a shop, i.e., the total time necessary for processing all jobs, under the assumption that we have control over the processing order of jobs, $v$, and the timing of events, $x_i(k)$, for all $i,k$.
In the remainder of this paper, we limit our focus to \textit{permutation flow shops}, where all jobs undergo the same sequence of processes in the same order and no job overtaking is allowed.
Mathematically, this corresponds to the following temporal ordering of events: $x_{i+1}(k)\geq x_i(k)$ and $x_i(k+1) \geq x_i(k)$ (i.e., $(A^0_{v_k})_{i+1,i}\geq 0$ and $(A^1_{v_k})_{i,i}\geq 0$) for all $i,k$.

In permutation flow shops described by SLDIs, the minimal makespan corresponds to the optimal value of the following optimization problem:
\begin{equation}\label{eq:problem_SLDIs}
    \begin{array}{cll}
        \min_{v\in S, x(1),\ldots,x(\quantities)\in\R^n} & & x_n(\quantities) - x_1(1) \\
        \mbox{subject to} & & \mbox{constraints \eqref{eq:SLDI}},
    \end{array}
\end{equation}
where the set $S$ is any subset of interest of the set of ordered sequences from $\Sigma$ of length $\quantities$.
For example, $S$ can be the set of permutations of $\Sigma$, in which case each job must be processed exactly once in the shop and $\quantities$ coincides with the number of elements in $\Sigma$.

Problem~\eqref{eq:problem_SLDIs} is, in general, NP-hard, and can be divided into two subproblems:
\begin{enumerate}[label=\textbf{P\arabic*}]
    \item\label{en:sp_makespan} the makespan computation, when $v$ is fixed;
    \item\label{en:sp_search} the search for the optimal $v\in S$ that attains the minimal makespan.
\end{enumerate}
In this contribution, we focus mainly on~\ref{en:sp_makespan}, which can be solved using linear programming, as all constraints~\eqref{eq:SLDI} are linear (see the last paragraph in the previous subsection), or graph-based approaches.
Here, we propose a technique based on max-plus algebra and inspired by the graph-theoretical interpretation of the problem. 

As for~\ref{en:sp_search}, various sophisticated strategies have been suggested for its solution in flow shops without upper bound constraints (where $B^0_{v_k}=B^1_{v_k} = \pazocal{T}$). 
Some of those can be adapted to the case study of the present paper; we mention branch and bound techniques for exact methods, heuristics and metaheuristics for approximations (see, e.g.,~\cite{FERNANDEZVIAGAS2017707}).
Clearly, the most trivial (and time-consuming) approach to solve~\ref{en:sp_search} in an exact way is an exhaustive search of all $|S|$ possible sequences $v\in S$.
In contrast, the case with upper bound constraints is less studied.
Some references are~\cite{fondrevelle2006permutation,samarghandi2019minimizing,chamnanlor2017embedding}; note that the class of shops described by~\eqref{eq:SLDI} is larger than the ones studied in these papers, in which $B^1_{v_k}$ is always $\pazocal{T}$.
This assumption is usually relaxed in hoist and cluster tool scheduling problems (e.g., see~\cite{KATS20081196,monch2011survey}).
Generally speaking, many manufacturing systems can be described by a sparse version of~\eqref{eq:SLDI}, i.e., where most of the elements in $A_{v_k}^0$ and $A_{v_k}^1$ are $-\infty$, and most of those in $B_{v_k}^0$ and $B_{v_k}^1$ are $+\infty$.

\subsection{Makespan computation}\label{su:makespan_computation}

\begin{figure*}
    \centering
    \resizebox{.87\textwidth}{!}{
        \begin{tikzpicture}[thick]

\newcommand{\rBread}{.12}
\newcommand{\wBread}{.8}
\newcommand{\hPipe}{.5}
\newcommand{\wMach}{2.5}
\newcommand{\hMach}{1}
\newcommand{\lMixerDivider}{.6}
\newcommand{\wMixer}{1.5}
\newcommand{\hMixer}{1.1}
\newcommand{\lMachines}{1.3}
\newcommand{\ifequals}[3]{\ifthenelse{\equal{#1}{#2}}{#3}{}}
\newcommand{\ccase}[2]{#1 #2} 
\newenvironment{switch}[1]{\renewcommand{\ccase}{\ifequals{#1}}}{}
\newcommand{\machName}[1]{\begin{switch}{#1}%
    \ccase{0}{Dividing}%
    \ccase{1}{Rounding}%
    \ccase{2}{Pre-proofing}%
    \ccase{3}{Rolling}%
\end{switch}}
\newcommand{\bread}[2]{
\draw [myyellow!70!black,fill=myyellow] (#1,#2) circle (\rBread);
}
\newcommand{\hOven}{3}
\newcommand{\wOven}{2}
\newcommand{\proofer}[2]{
    \draw [fill=red!10] (#1-\wOven/2,#2-\hOven/2) rectangle (#1+\wOven/2,#2+\hOven/2);

    \bread{#1-.2*\wOven}{#2-\hOven/4+\rBread}
    \bread{#1+.2*\wOven}{#2-\hOven/4+\rBread}
    \draw (#1-.4*\wOven,#2-\hOven/4) -- (#1+.4*\wOven,#2-\hOven/4);

    \bread{#1-.2*\wOven}{#2+\rBread}
    \bread{#1+.2*\wOven}{#2+\rBread}
    \draw (#1-.4*\wOven,#2) -- (#1+.4*\wOven,#2);

    \bread{#1-.2*\wOven}{#2+\hOven/4+\rBread}
    \bread{#1+.2*\wOven}{#2+\hOven/4+\rBread}
    \draw (#1-.4*\wOven,#2+\hOven/4) -- (#1+.4*\wOven,#2+\hOven/4);
}
\newcommand{\emptyProofer}[2]{
    \draw [fill=myblue!10] (#1-\wOven/2,#2-\hOven/2) rectangle (#1+\wOven/2,#2+\hOven/2);

    \draw (#1-.4*\wOven,#2-\hOven/4) -- (#1+.4*\wOven,#2-\hOven/4);

    \draw (#1-.4*\wOven,#2) -- (#1+.4*\wOven,#2);

    \draw (#1-.4*\wOven,#2+\hOven/4) -- (#1+.4*\wOven,#2+\hOven/4);
}
\newcommand{\oven}[2]{
    \draw [fill=red!10] (#1-\wOven/2,#2-\hOven/2) rectangle (#1+\wOven/2,#2+\hOven/2);

    \bread{#1-.2*\wOven}{#2-\hOven/4+\rBread}
    \bread{#1+.2*\wOven}{#2-\hOven/4+\rBread}
    \draw (#1-.4*\wOven,#2-\hOven/4) -- (#1+.4*\wOven,#2-\hOven/4);

    \bread{#1-.2*\wOven}{#2+\rBread}
    \bread{#1+.2*\wOven}{#2+\rBread}
    \draw (#1-.4*\wOven,#2) -- (#1+.4*\wOven,#2);

    \bread{#1-.2*\wOven}{#2+\hOven/4+\rBread}
    \bread{#1+.2*\wOven}{#2+\hOven/4+\rBread}
    \draw (#1-.4*\wOven,#2+\hOven/4) -- (#1+.4*\wOven,#2+\hOven/4);
}
\newcommand{\emptyOven}[2]{
    \draw [fill=red!10] (#1-\wOven/2,#2-\hOven/2) rectangle (#1+\wOven/2,#2+\hOven/2);

    \draw (#1-.4*\wOven,#2-\hOven/4) -- (#1+.4*\wOven,#2-\hOven/4);

    \draw (#1-.4*\wOven,#2) -- (#1+.4*\wOven,#2);

    \draw (#1-.4*\wOven,#2+\hOven/4) -- (#1+.4*\wOven,#2+\hOven/4);
}
%
\draw [myyellow!80!black,-stealth'] (0,\rBread+\hPipe+\hMixer+.5*\lMachines+.2) -- (0,\rBread+\hPipe+\hMixer+.2) node [black,pos=-.3] {Ingredients} node [black,pos=1.3] {Mixing};
\newcommand{\deltaMixer}{.6}
\fill [myyellow] (-\rBread,\rBread+\hPipe) -- ({-\rBread*(1-\deltaMixer)-\rBread*\deltaMixer-\wMixer/2*\deltaMixer},{(1-\deltaMixer)*\rBread+(1-\deltaMixer)*\hPipe+\rBread*\deltaMixer+\hPipe*\deltaMixer+\hMixer*\deltaMixer}) -- ({\rBread*(1-\deltaMixer)+\rBread*\deltaMixer+\wMixer/2*\deltaMixer},{\rBread*(1-\deltaMixer)+\hPipe*(1-\deltaMixer)+\rBread*\deltaMixer+\hPipe*\deltaMixer+\hMixer*\deltaMixer}) -- (\rBread,\rBread+\hPipe) -- cycle;

\draw (-\rBread,\rBread+\hPipe) -- (-\rBread-\wMixer/2,\rBread+\hPipe+\hMixer);
\draw (\rBread,\rBread+\hPipe) -- (\rBread+\wMixer/2,\rBread+\hPipe+\hMixer);
\draw [fill=myyellow] (0-\rBread,0+\rBread+\hPipe) -- (0-\rBread,0-\rBread) -- (\lMixerDivider,-\rBread) -- (\lMixerDivider,\rBread) -- (0+\rBread,0+\rBread) -- (0+\rBread,0+\rBread+\hPipe);

\foreach \i in {0,1,2,3}{

    \draw [fill=black!10] (\lMixerDivider+\i*\wMach+\i*\lMachines,-\hMach/2) rectangle (\lMixerDivider+\wMach+\i*\wMach+\i*\lMachines,\hMach/2) node [pos=.5] {\machName{\i}};
    \bread{\lMixerDivider+\wMach+\i*\wMach+\i*\lMachines+\lMachines/2}{-\hMach/2+\rBread}
    \draw (\lMixerDivider+\wMach+\i*\wMach+\i*\lMachines,-\hMach/2) -- (\lMixerDivider+\wMach+\lMachines+\i*\wMach+\i*\lMachines,-\hMach/2);
    \draw [myyellow!80!black,-stealth'] (\lMixerDivider+\wMach+\i*\wMach+\i*\lMachines+\lMachines/2-\lMachines/3,0) to [bend left] (\lMixerDivider+\wMach+\i*\wMach+\i*\lMachines+\lMachines/2+\lMachines/3,0);

}

\proofer{\lMixerDivider+\wMach+\lMachines+3*\wMach+3*\lMachines+\lMachines}{0}
\node at (\lMixerDivider+4*\wMach+5*\lMachines,1*\hOven/2+.4) {Proofing};

\draw [-stealth',myyellow!80!black] (\lMixerDivider+\wMach+\lMachines+3*\wMach+3*\lMachines+\lMachines+\wOven/2+\lMachines/6,0) to [bend left]  (\lMixerDivider+\wMach+\lMachines+3*\wMach+3*\lMachines+\lMachines+\wOven/2+5*\lMachines/6,0);

\oven{\lMixerDivider+4*\wMach+6*\lMachines+\wOven}{0}
\node at (\lMixerDivider+4*\wMach+6*\lMachines+\wOven,\hOven/2+.4) {Baking};

\draw [-stealth',myyellow!80!black] (\lMixerDivider+4*\wMach+6*\lMachines+\wOven+\wOven/2+.1*\lMachines,0) --  (\lMixerDivider+4*\wMach+6*\lMachines+\wOven+\wOven/2+.6*\lMachines,0) node [pos=1.2,black,align=center,anchor=west] {Finished\\products};

\end{tikzpicture}
    }
    \caption{
        Schematic representation of the bakery shop under study.
    }
    \label{fi:bakery_shop}
\end{figure*}
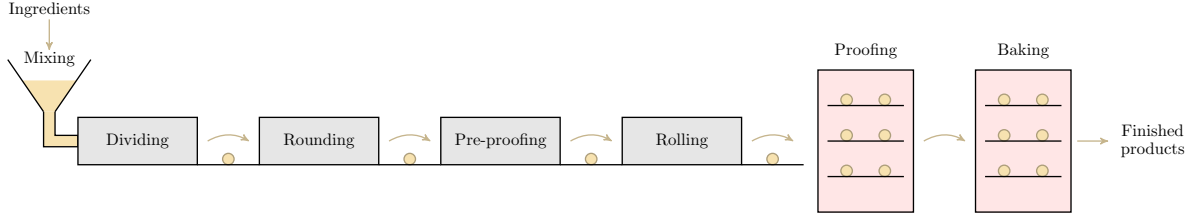

Let us consider Problem~\ref{en:sp_makespan}.
To solve it using max-plus algebra, we start rephrasing~\eqref{eq:SLDI} as a max-plus linear system of the form $A\otimes x \preceq x$ by means of the following proposition.

\begin{prop}[\cite{cuninghame2012minimax}]\label{pr:max-min}
Let $x,y\in\R^n$, $A,B\in\Rmax^{n\times n}$.
Then $x\preceq A^\sharp \stimes y \ \Leftrightarrow \ A\otimes x\preceq y$, and
\[
    \left\{
    \begin{array}{l}
        A\otimes x\preceq y\\
        B\otimes x\preceq y
    \end{array}
    \right.
    \quad
    \Leftrightarrow
    \quad
    (A\oplus B) \otimes x\preceq y.
\]
\end{prop}

Thanks to Proposition~\ref{pr:max-min}, it is possible to rewrite~\eqref{eq:SLDI} as 
\begin{equation}\label{eq:aux}
    \begin{array}{rc}
    \forall k\in\dint{1,\quantities}: & (A^0_{v_k}\oplus B^{0\sharp}_{v_k}) \otimes x(k) \preceq x(k),\\
    \forall k\in\dint{1,\quantities-1}: & A^1_{v_k} \otimes x(k) \preceq x(k+1),\\
    \forall k\in\dint{1,\quantities-1}: & B^{1\sharp}_{v_k} \otimes x(k+1) \preceq x(k).
    \end{array}
\end{equation}
We make the following substitutions: for all $k\in\dint{1,\quantities}$, $C_k \coloneqq A^0_{v_k}\oplus B^{0\sharp}_{v_k}$, for all $k\in\dint{1,\quantities-1}$, $I_k \coloneqq A^1_{v_k}$ and $P_k \coloneqq B^{1\sharp}_{v_k}$, $\tilde{x} \coloneqq [x^\top(1)\ \dots \ x^\top(\quantities)]^\top$, and
\begin{equation}\label{eq:Mv}
    M_v \coloneqq
    \begin{bmatrix}
        C_1 & P_1 & \pazocal{E} & \pazocal{E} & \cdots & \pazocal{E}\\
        I_1 & C_2 & P_2 & \pazocal{E} & \cdots & \pazocal{E}\\
        \pazocal{E} & I_2 & C_3 & P_3 & \cdots & \pazocal{E}\\
        \pazocal{E} & \pazocal{E} & I_3 & C_4 & \cdots & \pazocal{E}\\
        \vdots & \vdots & \vdots & \vdots & \ddots & \vdots\\
        \pazocal{E} & \pazocal{E} & \pazocal{E} & \pazocal{E} & \cdots & C_{\quantities}
    \end{bmatrix}\in \Rmax^{\quantities n\times \quantities n}.
\end{equation}
Then, \eqref{eq:aux} is equivalent to $M_v\otimes \tilde{x} \preceq \tilde{x}$, and the makespan can be computed from the next proposition.

\begin{prop}[\cite{butkovivc2010max}]
Given $A\in\Rmax^{n\times n}$, inequality $A\otimes x\preceq x$ admits a solution $x\in\R^n$ if and only if $\graph(A)\in\nonegset$. 
In this case (and only in this case), a solution of the optimization problem
\begin{equation}\label{eq:problem_max-plus}
    \begin{array}{cll}
        \min_{x\in\R^n} & & x_i-x_j \\
        \mbox{subject to} & & A\otimes x\preceq x.
    \end{array}
\end{equation}
exists and the minimal value of the cost function is $(A^*)_{ij}$.
\end{prop}
Moreover, it can be shown that vector $x = (A^*)_{\cdot j}$ (the $j$th column of $A^*$, for which $x_i = (A^*)_{ij}$ and $x_j = 0$) always solves~\eqref{eq:problem_max-plus}.
Therefore, if $\graph(M_v)\in\nonegset$, then the solution of Problem~\ref{en:sp_makespan} (i.e., the makespan) is equal to $(M_v^*)_{\quantities n,1}$, and $\tilde{x} = (M_v^*)_{\cdot 1}$ (the first column of $M_v^*$) corresponds to a consistent trajectory for which the makespan is attained. 

Even though this method solves our problem, it is not particularly fast: checking whether $\graph(M_v)\in\nonegset$ and computing $M_v^*$ take $\pazocal{O}((\quantities n)^3) = \pazocal{O}(\quantities^3 n^3)$ operations; when either $\quantities$ (i.e., the number of jobs to be processed) or $n$ (related to the number of processes in the shop) is large, this strategy becomes impractical. 
To speed up the computation, one must notice from~\eqref{eq:Mv} that matrix $M_v$ is sparse; only elements "close" to the diagonal of $M_v$ are not $-\infty$.
By exploiting this fact, the solution of Problem~\ref{en:sp_makespan} can be obtained in $\pazocal{O}(\quantities n^3)$, i.e., in linear time with respect to $\quantities$.
The algorithm comes directly from the following theorem, which is proven in the Appendix.

\begin{thm}\label{th:sparsity}
The solution of Problem~\ref{en:sp_makespan} exists if and only if, for all $k\in\dint{1,\quantities}$, $\graph(C_k)\in\nonegset$, and, for all $k\in\dint{1,\quantities -1}$, $\graph(\mathbb{C}_k)\in\nonegset$, where\footnote{Here and in the Appendix, we adopt the notation $A B = A\otimes B$.}
\[
    \mathbb{C}_i = \mathbb{P}_i(\mathbb{P}_{i+1}(\cdots(\mathbb{P}_{\quantities-1}\mathbb{I}_{\quantities-1})^*\cdots)^*\mathbb{I}_{i+1})^*\mathbb{I}_i,
\]
$\mathbb{P}_i = C_i^* P_i C_{i+1}^*$, and $\mathbb{I}_i = C_{i+1}^*I_i C_i^*$.
If the conditions above hold, then the makespan is given by element $(n,1)$ of matrix
\[
    \mathbb{M} = \mathbb{I}_{\quantities-1} \mathbb{C}_{\quantities-1}^* \mathbb{I}_{\quantities-2} \mathbb{C}_{\quantities-2}^* \cdots \mathbb{I}_1\mathbb{C}_1^*.
\]
\end{thm}

\begin{rem}\label{re:sparsity}
In flow shops, matrices $P_k$, $I_k$, and $C_k$ are typically sparse as well.
This observation, partially exploited in the case study in the next section, leads to an even faster makespan computation.
\end{rem}

\begin{rem}
The complexity for solving~\eqref{eq:problem_SLDIs} using the makespan computation proposed in the present section and an exhaustive search of $v$ is $\pazocal{O}(|S|\quantities n^3)$.
\end{rem}

\section{Case study}\label{se:4}

We start this section by giving a mathematical description, in standard algebra, of the industrial bakery shop under study, which is schematically represented in Figure~\ref{fi:bakery_shop}.
After that, the problem of minimizing its makespan is considered and solved in the max-plus algebra.
Finally, comparisons with other algorithms are made.

\subsection{The bakery shop}\label{su:bakery_shop}

It consists of $\machines=7$ stages in series, each containing one machine: 
(1) mixing of ingredients, (2) dough dividing, (3) rounding, (4) pre-proofing, (5) rolling, (6) proofing, and (7) baking.
The bakery shop can process $\types$ types of products (or jobs).
The quantity of products of type $\type\in\dint{1,\types}$ to be processed is denoted by $\quantities_\type\in\No$; this value may change every day, according to customers' demand.
The total quantity to be processed in a day is $\quantities \coloneqq \sum_{\type=1}^{\types} \quantities_{\type}$.

Each product undergoes the same sequence of processes in the same order, as indicated in Figure~\ref{fi:bakery_shop}, but the processing time in a machine is different for different product types; a manufacturing system of this kind is called flow shop (\cite{pinedo2016scheduling}).
We number products progressively from $k=1$ to $k=\quantities$ according to their entrance time into the bakery shop, and we denote by $\type(k)\in\dint{1,\types}$ the type of the $k$th product.
Let $\xi_{\machine}(k)\in\R_{\geq 0}$ and $\xi_{\machine}'(k)\in\R_{\geq \xi_{\machine}(k)}$ indicate, respectively, the time instants at which the $k$th product enters and leaves machine $\machine\in\machineset$; these times are collected in the vector $x(k) \coloneqq [\xi_1(k)\ \xi_1'(k)\ \ldots \ \xi_{\machines}(k)\ \xi_{\machines}'(k)]^\top\in\R^{2\machines}_{\geq 0}$. 
We denote by $\tau_{\machine\type}^-$ and $\tau_{\machine\type}^+$ the minimum and maximum time for processing a product of type $\type\in\typeset$ in machine $\machine\in\machineset$, respectively, and by $\tau_{\machine}^-$ and $\tau_{\machine}^+$ the minimum and maximum time to transport products from machine $\machine$ to $\machine+1$, with $\machine\in\dint{1,\machines-1}$. 
Then, the following inequalities must hold for all $k\in\quantityset$:
\begin{equation}
\label{eq:1}
\begin{array}{ccr}
     \tau_{\machine\type(k)}^- \leq \xi_\machine'(k) - \xi_{\machine}(k) \leq \tau_{\machine\type(k)}^+ & & \forall \machine\in\machineset, 
     \\
    \tau_{\machine}^- \leq \xi_{\machine+1}(k) - \xi_{\machine}'(k) \leq \tau_{\machine}^+ & & \forall \machine\in\dint{1,\machines-1}. 
\end{array}
\end{equation}
When not otherwise stated, it is assumed that $0< \tau_{\machine\type}^-<\tau_{\machine\type}^+<+\infty$ and $0< \tau_{\machine}^-<\tau_{\machine}^+<+\infty$ for all $\machine,\type$.

The dynamics of products in the shop varies qualitatively from machine to machine and, based on their common characteristics, machines can be divided in three groups.
In the following we describe the stages of the shop -- for convenience, we proceed in reverse order.

\textbf{Proofing and baking.}
The proofer and the oven have the same capacity, i.e., they can process the same number of products at a time; however, their capacity, denoted by $\capacity_\type\in\N$, varies with the type $\type$ of product to be processed, since the larger the surface area of a product, the fewer products can fit in a machine.

Products of different type require different temperatures during proofing and baking; thus, it is assumed that only products of the same type can reside in a proofer or oven at the same time.
Moreover, once the proofing or baking process starts, products cannot be inserted to or removed from the machine until the process is completed; therefore, products must enter or leave these machines in batches.
The number of batches of products of type $\type$ that will enter the proofer and oven on a given day is $\batches_{\type}\coloneqq \ceil{\frac{\quantities_{\type}}{\capacity_{\type}}}$, and the total number of batches to be processed in the shop is $\batches\coloneqq \sum_{\type=1}^{\types} \batches_{\type}$.
The first $\floor{\frac{\quantities_\type}{\capacity_\type}}$ batches of products of type $\type$ consist of $\capacity_\type$ products (i.e., they fill up the machines completely), whereas the last batch consists of $\quantities_\type \mod \capacity_\type$ products (i.e., if $\quantities_\type$ is not divisible by $\capacity_\type$, the last batch does not fill up the machines completely).
In the following, the term "batch" will be used to indicate the quantity of ingredients or dough corresponding to either $\quantities_\type \mod \capacity_\type$ or $\capacity_\type$ products of type $\type$, depending on whether we are referring to the last batch of type $\type$ or not.

We let the batch of type $\type(k)$ to which the $k$th product belongs be denoted by $\batch(k) \in\dint{1,\batches_{\type(k)}}$.
Thus, products $k$ and $k'$ belong to a different batch if either $\batch(k)\neq \batch(k')$ or $\type(k)\neq \type(k')$.
According to the description above, the dynamics of products in the proofing and baking stages can be formulated mathematically as follows: for all $k\in\dint{1,\quantities-1}$, if $\batch(k+1)=\batch(k)$ and $\type(k+1) = \type(k)$, then\footnote{Note that $0 \leq b - a \leq 0$ is equivalent to $a = b$.
Moreover, observe that here upper-bound constraints are not used to represent quality specifications, but rather a batching phenomenon.}
\begin{equation}
    \begin{array}{llr}
0 \leq \xi_\machine(k+1) - \xi_\machine(k) \leq 0 && \forall \machine\in\dint{6,7},\\
0\leq \xi_\machine'(k+1) - \xi_\machine'(k) \leq 0 && \forall \machine\in\dint{6,7};
    \end{array}
\end{equation}
otherwise, 
\begin{align}
0 \leq \xi_\machine(k+1) - \xi_\machine'(k) && \forall \machine\in\dint{6,7}.
\end{align}
The latter inequality forces batches to enter the proofer or the oven only after the removal of the previous batch.

From the rolling to the proofing stage and from the proofing to the baking stage, products are transported in batches by means of trolleys.
It is assumed that the number of trolleys in the shop is sufficiently high, so that there is always one available when needed.

\textbf{Dividing, rounding, pre-proofing, and rolling.}
Each machine in these stages has unitary capacity; thus, for all $k\in\dint{1,\quantities-1}$,
\begin{align}
     0 \leq \xi_\machine(k+1) - \xi_\machine'(k) && \forall \machine\in\dint{2,5}.
\end{align}
Between these machines there is no intermediate storage, i.e., when a process terminates the next one starts; hence, $\tau_{\machine}^-=\tau_{\machine}^+ = 0$ for all $\machine\in\dint{2,4}$. 
Moreover, due to technological constraints, products are not allowed to wait in these machines after being processed\footnote{In scheduling theory, this is referred to as \textit{no-wait} requirement.}, i.e., $\tau_{\machine\type}^{-}=\tau_{\machine\type}^+$ for all $\machine\in\dint{2,5}$, $\type\in\typeset$.

\textbf{Mixing.}
The capacity of the mixer is assumed to be sufficiently large to contain an arbitrary quantity of ingredients. 
Ingredients in quantity corresponding to a batch of products are inserted all at once in the mixer, and between the insertion of products of different types, the mixer needs to be cleaned for a time equal to $\tau_{\textup{clean}}$.
Formally, for all $k\in\dint{1,\quantities-1}$, if $b(k+1) = b(k)$ and $\type(k+1)=\type(k)$,
\begin{align}
    0 \leq \xi_1(k+1) - \xi_1(k) \leq 0,
\end{align}
otherwise, if $\type(k+1) \neq \type(k)$,
\begin{align}\label{eq:cleaning}
    \tau_{\textup{clean}} \leq \xi_1(k+1) - \xi_1'(k).
\end{align}
(In the case when $\batch(k+1)\neq \batch(k)$ and $\type(k+1)=\type(k)$, no additional constraint is required.)
Without loss of generality, we assume that the mixer operates in a FIFO manner, i.e., for all $k<k'$,
\begin{align}\label{eq:7}
    0 \leq \xi_1'(k') - \xi_1'(k).
\end{align}

The dividing machine, used to divide the dough into smaller pieces of the size of a product, is physically connected to the mixer; thus, $\tau_{1}^-=\tau_{1}^+=0$.

\begin{rem}
It is easy to verify that~(\ref{eq:1}--\ref{eq:7}) and the definition of $k$ imply that the bakery shop is a permutation flow shop (see Section~\ref{su:makespan_permutation}).
\end{rem}

\subsection{The scheduling problem in the bakery shop}

The quantity of products of each type to be processed in the bakery shop is decided every day according to the market demand.
Our goal is to determine, for each product, the optimal start and end times of each process, i.e., to minimize the makespan.
The difficulty of this problem lies in the fact that, to be able to use the optimal production plan in the shop, a solution must be found before the work shift starts; this imposes a maximum acceptable time for computing the solution of about 15 minutes.

Before approaching the problem, we make the following reasonable assumption, which narrows the search space:
after starting processing a product of type $\type$ in the shop, all $\quantities_\type$ products of that type need to be processed before switching to a different type.
The assumption is justified by the fact that, after switching to a different product type, it is necessary to perform a time-consuming cleaning procedure of the mixer (see~\eqref{eq:cleaning}). 
Therefore, only production plans where the number of cleaning procedures is minimal are considered.

\begin{figure}
    \resizebox{1\linewidth}{!}{
    \input{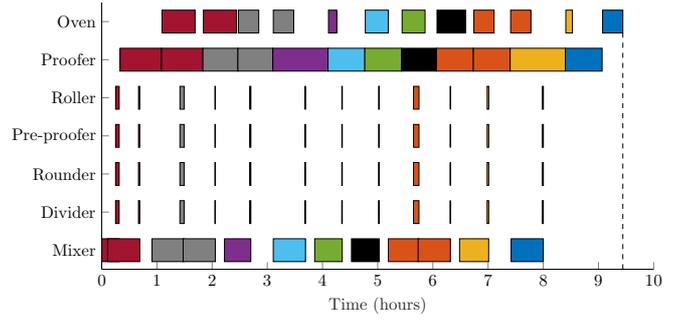}
    }
    \caption{Gantt chart of a consistent trajectory for the bakery shop. Each rectangle corresponds to a batch of products being processed in a machine, different colors indicate different product types, and the dashed line shows the makespan ($9.4$ hours); in this example $\types=9$ and $\batches=12$.}
    \label{fi:gantt}
\end{figure}

Let a schedule $w = (\type_1,\ldots,\type_{\types})$ indicate the order in which product types enter the first stage of the bakery, i.e., the mixing stage.
Due to the assumption above, $w$ is a permutation of set $\typeset$ and, once a schedule $w$ is fixed, the type $\type(k)$ and batch $\batch(k)$ of product $k\in\quantityset$ are uniquely determined: given $i\in\typeset$,
\begin{equation}\label{eq:type}
    \type(k) = \type_{i} \mbox{ for all } k\in\left\llbracket\sum_{h=1}^{i-1}\quantities_{\type_h}+1,\sum_{h=1}^{i}\quantities_{\type_h}\right\rrbracket,
\end{equation}
\begin{equation}\label{eq:assignment}
    \footnotesize
    \batch(k) = 
    \begin{dcases}
        1 & \mbox{for } k\in\dint{\sum_{h=1}^{i-1}\quantities_{\type_h}+1,\sum_{h=1}^{i-1}\quantities_{\type_h}+\capacity_{\type_i}}\\
        \cdots & \\
        \batches_{\type_i} & \mbox{for } k\in\dint{\sum_{h=1}^{i-1}\quantities_{\type_h}+(\batches_{\type_i}-1)\capacity_{\type_i}+1,\sum_{h=1}^{i}\quantities_{\type_h}}.
    \end{dcases}
\end{equation}
Figure~\ref{fi:gantt} shows an example of consistent trajectory for the bakery shop.
Note that batches of the same product type are processed consecutively, in accordance with the above assumption.

The scheduling problem is stated as follows:
\begin{equation}\label{eq:problem}
    \begin{array}{cll}
        \min_{w,x(1),\dots,x(\quantities)} & & \xi_7'(\quantities)-\xi_1(1) \\
        \mbox{subject to} & & \mbox{constraints (\ref{eq:1}--\ref{eq:7})},\\
        & & \mbox{assignments (\ref{eq:type}--\ref{eq:assignment})}.
    \end{array}
\end{equation}
Quantity $\xi_7'(\quantities)-\xi_1(1)$ corresponds to the time difference between the exit of the last batch from the baking stage (i.e., the last stage of the shop) and the shop start-up; thus, it coincides with the shop's makespan.

\subsection{Solution in the max-plus algebra}\label{su:solution_bakery}

We start by rewriting the dynamics of products in the bakery shop as SLDIs.
Consider the set of modes $\Sigma = \bigcup_{\type\in\typeset}\{\wA_{\type},\wB_{\type},\wC_{\type}\}$.
Given a schedule $w = (\type_1,\ldots,\type_{\types})$, we build a sequence of modes $v$ of length $\quantities$ such that:
\begin{itemize}
    \item $v_k\in\{\wA_\type,\wB_\type,\wC_\type\}$ if the type of the $k$th product is $\type$ ($\type(k)=\type$),
    \item $v_k = \wA_\type$ if the $k$th and the $(k+1)$st products belong to the same batch ($\batch(k+1)=\batch(k)$ and $\type(k+1) = \type(k)$),
    \item $v_k = \wB_\type$ if the $k$th and the $(k+1)$st products belong to different batches of the same type ($\batch(k+1)\neq\batch(k)$ and $\type(k+1) = \type(k)$),
    \item $v_k = \wC_\type$ if the $k$th and the $(k+1)$st products are of different types ($\type(k+1) \neq \type(k)$).
\end{itemize}
It is then easy to rewrite inequalities (\ref{eq:1}--\ref{eq:7}) in the form of SLDIs by constructing matrices $A^0, A^1, B^0, B^1$ for each mode in $\Sigma$.

By applying directly Theorem~\ref{th:sparsity}, combined with an exhaustive search of schedule $w$, we can finally compute the minimal makespan of the bakery shop in time $\pazocal{O}(\types! \quantities \machines^3)$ (as the number of possible $w$ is $\types !$ and the dimension of vector $x(k)$ is $2\machines$).
Actually, the complexity can be further reduced in the considered case study: a careful analysis of (\ref{eq:1}--\ref{eq:7}) reveals that $B^1_{\wC_\type}=\pazocal{T}$ for all $\type\in\typeset$. 
The consequence is that formulas from Theorem~\ref{th:sparsity} simplify, as some of the terms cancel out; in combination with the simple assignments~(\ref{eq:type}--\ref{eq:assignment}), this allows to pre-compute many of the terms composing matrix $\mathbb{M}$ in Theorem~\ref{th:sparsity}.
Overall, time complexity reduces to $\pazocal{O}((\quantities + \types ! \types) \machines^3)$, which approximates $\pazocal{O}(\types ! \types \machines^3)$ when $\quantities \ll \types ! \types$; hence, the contribution of the quantity of products to be produced, $\quantities$, on the computational time becomes negligible when the number of product types, $\types$, is sufficiently high.

\subsection{Comparison of different techniques}

We tested four algorithms for the makespan computation in the bakery shop on a PC with an Intel i7 processor at 2.20Ghz.
The considered scenario corresponds to a typical daily production plan in the bakery: a total number of $\quantities=975$ products of $\types=9$ different types are processed in $\batches=12$ batches (as in the case of Figure~\ref{fi:gantt}).
The comparison is made between a linear programming solver (CPLEX's dual-simplex method), a graph-based method (Bellman-Ford shortest path algorithm), the algorithm derived by the direct application of Theorem~\ref{th:sparsity} (of complexity $\pazocal{O}(\types! \quantities \machines^3)$), and the faster approach described in Section~\ref{su:solution_bakery} (of complexity approximately $\pazocal{O}(\types ! \types \machines^3)$).
For the makespan computation, the first two methods have been implemented to solve Problem~\eqref{eq:problem_max-plus}; whenever possible, the sparsity of the problem has been taken into account during implementation.

The following results have been obtained: the computation of the makespan for a single schedule takes $7.25\cdot 10^{-2}$s using the dual-simplex method, $2.99\cdot 10^{-2}$s using the Bellman-Ford algorithm, $1.80\cdot 10^{-2}$s using Theorem~\ref{th:sparsity}, and $3.77\cdot 10^{-5}$s using the method from Section~\ref{su:solution_bakery}.
The latter approach thus saves 3 orders of magnitude in terms of time.
With an exhaustive search, the optimal makespan is found, respectively, in $7.30$ hours, $3.02$ hours, $1.82$ hours, and $13.7$ seconds.
Hence, only the latter approach meets the time limitation of 15 minutes for solving the problem.

When more than $\types=10$ types of product are considered, the exhaustive search becomes too slow even for the faster approach from Section~\ref{su:solution_bakery}; indeed, for $\types=11$ and the same values of $\quantities$ and $\batches$ as before, $25.1$ minutes are needed to compute the makespan for all the schedules.
Thus, for $\types>10$, combining the formulas from Theorem~\ref{th:sparsity} with branch and bound methods, heuristics or meta-heuristics becomes necessary.

\section{Conclusions and future work}\label{se:conclusions}

In this contribution, we modeled a section of an industrial bakery as SLDIs and, based on the max-plus algebra, we obtained a closed formula for the makespan in flow shops with time-window constraints.
The formula has the advantage of revealing sparsity patterns otherwise obscured by the black-box nature of optimization problems; in the present case study, hours of computations are avoided thanks to it.

In practice, bakery shops are often hybrid flow shops, i.e., each stage of the shop may contain several machines in parallel.
Scheduling problems in hybrid flow shops are significantly more challenging compared to flow shops; in future work, we aim to extend our analysis to this scenario.
\section*{Appendix}

In this section, we prove Theorem~\ref{th:sparsity}.
All symbols used in the following propositions are as in Section~\ref{su:makespan_computation}.
We start from the following lemma.

\begin{lem}\label{le:block_star}
Let $a\in\Rmax^{n_1\times n_1}$, $b\in\Rmax^{n_1\times n_2}$, $c\in\Rmax^{n_2\times n_1}$, and $d\in\Rmax^{n_2\times n_2}$.
Then
\[
    \begin{bmatrix}
        a & b\\c & d
    \end{bmatrix}^* = 
    \begin{bmatrix}
        a^*(a^*bd^*ca^*)^*a^* & a^*bd^*(d^*ca^*bd^*)^*\\
        d^*ca^*(a^*bd^*ca^*)^* & d^*(d^*ca^*bd^*)^*d^*
    \end{bmatrix}.
\]
\end{lem}
\begin{pf}
    The lemma can be proven by applying formulas (6.3) -- (6.10) from~\cite{hardouin2018control} to Algorithm 2 from the same paper for the computation of the Kleene star.
    \qed
\end{pf}

In the following, we will use the notation
\[
    M_v^* = 
    \begin{bsmallmatrix}
        \mathcal{M}^{11} & \cdots & \mathcal{M}^{1\quantities}\\
        \svdots & \sddots & \svdots \\
        \mathcal{M}^{\quantities 1} & \cdots & \mathcal{M}^{\quantities \quantities}
    \end{bsmallmatrix},
\]
where each $\mathcal{M}^{ij}$ is a matrix of dimension $n\times n$, and $\mathcal{M}^{\quantities 1}$ coincides with $\mathbb{M}$ from Theorem~\ref{th:sparsity}.

\begin{prop}\label{pr:aux_C1}
    $\mathcal{M}^{11} = C_1^* \mathbb{C}_1^* C_1^*$.
\end{prop}
\begin{pf}
The proof is done by induction on $\quantities$.
For $\quantities = 1$ and $\quantities = 2$, the formula comes from the previous lemma.
Suppose that the equality holds for $\quantities-1$, with $\quantities\geq 3$; we prove that it holds also for $\quantities$.
By partitioning matrix $M_v$ in four blocks, such that the upper-left block coincides with $C_1$, from the previous lemma we get the following expression for $\mathcal{M}^{11}$:
\[
    C_1^*\left(C_1^*\begin{bmatrix}P_1 & \pazocal{E} & \cdots & \pazocal{E}\end{bmatrix}
    \begin{bmatrix}
        C_2 & P_2 & \cdots & \pazocal{E}\\
        I_2 & C_3 & \cdots & \pazocal{E}\\
        \vdots & \vdots & \ddots & \vdots\\
        \pazocal{E} & \pazocal{E} & \cdots & C_{\quantities}
    \end{bmatrix}^*
    \begin{bmatrix}I_1 \\ \pazocal{E} \\ \vdots \\ \pazocal{E}\end{bmatrix}
    C_1^*\right)^*C_1^*.
\]
By the induction hypothesis, the latter formula simplifies to $C_1^*(C_1^* P_1 C_2^* \mathbb{C}_2^* C_2^* I_1 C_1^*)^*C_1^*$, which coincides with $C_1^*\mathbb{C}_1^*C_1^*$.
\qed
\end{pf}

\begin{prop}\label{pr:aux_MK1}
    $\mathcal{M}^{\quantities 1} = \mathbb{I}_{\quantities-1}\mathbb{C}_{\quantities -1}^* \mathbb{I}_{\quantities -2}\mathbb{C}_{\quantities-2}^* \cdots \mathbb{I}_1 \mathbb{C}_{1}^*$.
\end{prop}
\begin{pf}
The proof is similar to the previous one.
For $\quantities=1$ or $\quantities=2$, it comes from Lemma~\ref{le:block_star}.
Assuming the formula is correct for $\quantities -1$ with $\quantities \geq 3$, we can prove its correctness for $\quantities$ as follows.
We partition $M_v$ as in the previous proof; then, $\mathcal{M}^{\quantities 1}$ is equal to the bottom matrix block of
\[
    \begin{bmatrix}
        C_2 & P_2 & \cdots & \pazocal{E}\\
        I_2 & C_3 & \cdots & \pazocal{E}\\
        \vdots & \vdots & \ddots & \vdots\\
        \pazocal{E} & \pazocal{E} & \cdots & C_{\quantities}
    \end{bmatrix}^*
    \begin{bmatrix}
        I_1 \\ \pazocal{E} \\ \vdots \\ \pazocal{E}
    \end{bmatrix}
    C_1^*
    \Biggl(C_1^* \begin{bmatrix}P_1 & \pazocal{E} & \cdots & \pazocal{E}\end{bmatrix} \otimes
\]\[
    \otimes \begin{bmatrix}
        C_2 & P_2 & \cdots & \pazocal{E}\\
        I_2 & C_3 & \cdots & \pazocal{E}\\
        \vdots & \vdots & \ddots & \vdots\\
        \pazocal{E} & \pazocal{E} & \cdots & C_{\quantities}
    \end{bmatrix}^*
    \begin{bmatrix}
        I_1 \\ \pazocal{E} \\ \vdots \\ \pazocal{E}
    \end{bmatrix}
    C_1^*\Biggr)^*.
\]
Expanding the expression using the induction hypothesis and Proposition~\ref{pr:aux_C1}, we get the desired formula. 
\qed
\end{pf}

\begin{figure}
    \centering
    \resizebox{.7\linewidth}{!}{
    \begin{tikzpicture}

\node [place] (n1) {$\overline{1}$};
\node [place,right=of n1] (n2) {$\overline{2}$};
\node [place,right=of n2] (n3) {$\overline{3}$};
\node [place,right=of n3] (n4) {$\overline{4}$};
\node [place,right=of n4] (n5) {$\overline{5}$};

\draw [-stealth'] (n1) to [out=90-20,in=90+20,loop] node [above] {$C_1$} (n1);
\draw [-stealth'] (n2) to [out=90-20,in=90+20,loop] node [above] {$C_2$} (n2);
\draw [-stealth'] (n3) to [out=90-20,in=90+20,loop] node [above] {$C_3$} (n3);
\draw [-stealth'] (n4) to [out=90-20,in=90+20,loop] node [above] {$C_4$} (n4);
\draw [-stealth'] (n5) to [out=90-20,in=90+20,loop] node [above] {$C_5$} (n5);

\draw [-stealth'] (n1) to [bend left] node [above] {$I_1$} (n2);
\draw [-stealth'] (n2) to [bend left] node [above] {$I_2$} (n3);
\draw [-stealth'] (n3) to [bend left] node [above] {$I_3$} (n4);
\draw [-stealth'] (n4) to [bend left] node [above] {$I_4$} (n5);

\draw [-stealth'] (n2) to [bend left] node [below] {$P_1$} (n1);
\draw [-stealth'] (n3) to [bend left] node [below] {$P_2$} (n2);
\draw [-stealth'] (n4) to [bend left] node [below] {$P_3$} (n3);
\draw [-stealth'] (n5) to [bend left] node [below] {$P_4$} (n4);

\end{tikzpicture}}
    \caption{Lumped-node representation of $\graph(M_v)$\\for $\quantities=5$.}
    \label{fi:graph}
\end{figure}

\begin{prop}\label{pr:nonegset}
$\graph(M_v)\in\nonegset$ if and only if, for all $k\in\dint{1,\quantities}$, $\graph(C_k)\in\nonegset$ and, for all $k\in\dint{1,\quantities -1}$, $\graph(\mathbb{C}_k)\in\nonegset$.
\end{prop}
\begin{pf}
We recall that $\graph(M_v) \in\nonegset$ iff there are no elements equal to $+\infty$ in the diagonal of $M_v^*$.
Before continuing, it is useful to visualize the structure of $\graph(M_v)$; the lumped-node representation of Figure~\ref{fi:graph} illustrates it in the case $\quantities = 5$.
In this simplified representation, $\overline{k}$ indicates the set of nodes $\{(k-1)n+1,\ldots,kn\}$ for any $k\in\dint{1,5}$, and the matrix, say $Y$, associated to an arc from $\overline{k}_1$ to $\overline{k}_2$ indicates that an arc in $\graph(M_v)$ from node $(k_1-1)n+j$ to node $(k_2-1)n+i$ exists iff $Y_{ij}\neq -\infty$, and that its weight is equal to $Y_{ij}$, for all $i,j\in\dint{1,n}$.

We partition the set of circuits of $\graph(M_v)$ in $\quantities$ subsets as follows: (1) circuits visiting at least one node in $\overline{1}$, (2) circuits that do not visit any node in $\overline{1}$, but that visit at least one node in $\overline{2}$, \dots, ($\quantities$) circuits that do not visit any node in $\overline{\quantities -1}$, but that visit at least one node in $\overline{\quantities}$.
Let $k\in\dint{1,\quantities-1}$.
From Proposition~\ref{pr:aux_C1}, the maximal weights of circuits in the $k$th subset correspond to elements in the diagonal of $C_k^*\mathbb{C}_k^* C_k^*$; thus, these are non-positive iff $\graph(C_k)\in\nonegset$ and $\graph(\mathbb{C}_k)\in\nonegset$.
As for circuits in the $\quantities$th subset, their maximal weights come from the diagonal elements of $C_{\quantities}^*$, which are finite iff $\graph(C_{\quantities})\in\nonegset$.
\qed
\end{pf}
Theorem~\ref{th:sparsity} is an immediate consequence of Propositions~\ref{pr:aux_MK1} and~\ref{pr:nonegset}.

\bibliography{references}

\end{document}